\documentclass[12pt,dvipdfmx]{article}
\usepackage{tikz}
\usetikzlibrary{intersections,calc,arrows}
\usepackage{amssymb}
\date{}
\makeatletter
\newcommand{\figcaption}[1]{\def\@captype{figure}\caption{#1}}
\newcommand{\tblcaption}[1]{\def\@captype{table}\caption{#1}}

\@addtoreset{equation}{section}
\makeatother
\newcommand{\qed}{\hbox{\rule[-2pt]{3pt}{6pt}}}

\begin{document}
\title {\bf Asymptotic behavior of bifurcation curves of one-dimensional nonlocal elliptic equations}

\author{Tetsutaro Shibata 
\\
Laboratory of Mathematics, School of Engineering
\\
Graduate School of Advanced Science and Engineering
\\
Hiroshima University, 
Higashi-Hiroshima, 739-8527, Japan
}

\maketitle
\footnote[0]{E-mail: tshibata@hiroshima-u.ac.jp}
\footnote[0]{This work was supported by JSPS KAKENHI Grant Number JP21K03310.}

\begin{abstract}
We study the one-dimensional nonlocal elliptic equation
\begin{eqnarray*}
-\left(\int_0^1 \vert u(x)\vert^p dx + b\right)^q u''(x)&=& \lambda u(x)^p , 
\enskip x \in I:= (0,1), \enskip u(x) > 0, 
\enskip x\in I, 
\\
u(0) &=& u(1) = 0,
\end{eqnarray*}
where $b \ge 0, p \ge 1, q > 1 - \frac{1}{p}$ are given constants and $\lambda > 0$ is a bifurcation parameter.  
We establish the global behavior of bifurcation diagrams and 
precise asymptotic formulas for $u_\lambda(x)$ as $\lambda \to \infty$.

\end{abstract}



\section{Introduction} 		      

We consider the following one-dimensional nonlocal elliptic equation

\begin{equation}
\left\{
\begin{array}{l}
-\left(\displaystyle{\int_0^1} \vert u(x)\vert^p dx + b\right)^q u''(x)= \lambda u(x)^p, \enskip 
x \in I:= (0,1),
\vspace{0.1cm}
\\
u(x) > 0, \enskip x\in I, 
\vspace{0.1cm}
\\
u(0) = u(1) = 0,
\end{array}
\right.
\end{equation}
where $b, p, q$ are given constants satisfying 

\begin{eqnarray}
b \ge 0, \quad p \ge 1, \quad q > 1 - \frac{1}{p}
\end{eqnarray}
and $\lambda > 0$ is a bifurcation parameter. 
(1.1) is the model equation of the following nonlocal problem (1.3) considered in [7]. 
\begin{equation}
\left\{
\begin{array}{l}
-a\left(\displaystyle{\int_0^1} \vert u(x)\vert^p dx\right)u''(x)= h(x) f(x,u(x)), \enskip 
x \in I,
\vspace{0.1cm}
\\
u(x) > 0, \enskip x\in I, 
\vspace{0.1cm}
\\
u(0) = u(1) = 0,
\end{array}
\right.
\end{equation}
where $a= a(w)$ is a real-valued continuous function. Let 
$\Vert u\Vert_p := \left(\int_I \vert u(x)\vert^p dx\right)^{1/p}$. If we put 
$a(\Vert u\Vert_p) = (\Vert u\Vert_p^p + b)^q$, $h(x) \equiv \lambda$ and $f(x,u) = u^p$ 
in (1.3), then we obtain (1.1). 
Nonlocal elliptic problems as (1.3) have been studied intensively by many authors, 
since they arise in various physical models. We refer to [4,8,9,11] and the refernces therein. As for standard and nonlocal bifurcation problems, there are many results for the equations with 
different types of nonlinearities. We refer to [2, 5] and the references therein. The purpose of this paper is to obtain the 
precise asymptotic behavior of bifurcation curves $\lambda = \lambda(\alpha)$ and $u_\lambda$ as $\lambda \to \infty$ by 
focusing on the typical nonlocal problem (1.1). Here, 
$\alpha:= \alpha_\lambda =\Vert u_\lambda\Vert_\infty$ for given $\lambda > 0$.

To state our results, we prepare the following notation. For $p > 1$, let 
\begin{equation}
\left\{
\begin{array}{l}
-W''(x)= W(x)^p, \enskip 
x \in I,
\vspace{0.1cm}
\\
W(x) > 0, \enskip x\in I, 
\vspace{0.1cm}
\\
W(0) = W(1) = 0.
\end{array}
\right.
\end{equation}
We know from [6] that there exists a unique solution $W_p(x)$ of (1.4).

\vspace{0.2cm}

\noindent
{\bf Theorem 1.1.} {\it Let $b = 0$ in (1.1). Then there exists a unique solution 
$u_\lambda$ of (1.1) for any given $\lambda > 0$. Furthermore, the following 
formulas hold. 

\noindent
(i) Assume that $p > 1$. Then for a given $\lambda > 0$,
\begin{eqnarray}
\lambda &=& 2(p+1)L_p^{2-q}M_p^q\alpha^{pq-p+1},
\\
u_\lambda(x)&=& \lambda^{1/(pq-p+1)}\left\{(2(p+1))^{p/(p-1)}L_p^{(p+1)/(p-1)}
M_p\right\}^{-q/(pq-p+1)}W_p(x),
\end{eqnarray}
where 
\begin{eqnarray}
L_p:= \int_0^1 \frac{s}{\sqrt{1-s^{p+1}}}ds, \enskip
 M_p:= \int_0^1 \frac{s^p}{\sqrt{1-s^{p+1}}}ds.
 \end{eqnarray}
(ii) Assume that $p = 1$. Then 
\begin{eqnarray}
\lambda = 2^q\pi^{2-q}\alpha^q.
\end{eqnarray}
}

We next consider the case $b > 0$. To clarify our intention, we start from the 
simplest case $p = 2$ and $q = 1$. We have
\begin{eqnarray}
\Vert W_p\Vert_p &=& (2(p+1))^{1/(p-1)}M_p^{1/p}L_p^{(p+1)/(p(p-1))} \qquad (p > 1),
\\
\Vert W_p\Vert_\infty &=&  (2(p+1))^{1/(p-1)}L_p^{2/(p-1)}.
\end{eqnarray} 
We will show (1.9) and (1.10) in the last part of Section 2. 

\vspace{0.2cm}

\noindent
{\bf Theorem 1.2.} {\it Let $b > 0$, $p = 2, q = 1$ and 
\begin{eqnarray}
\lambda_0:= 2b^{1/2}\Vert W_2\Vert_2.
\end{eqnarray} 
(i) If $0 < \lambda < \lambda_0$, then there exists no solution of (1.1).

\noindent
(ii) If $\lambda = \lambda_0$, then (1.1) has a unique solution 
\begin{eqnarray}
u_\lambda(x) = \frac{\lambda_0}{2}\Vert W_2\Vert_2^{-2}W_2(x).
\end{eqnarray}

\noindent
(iii) If $\lambda > \lambda_0$, then there exist exactly two solutions 
$u_{1,\lambda}, u_{2,\lambda}$ of (1.1) such that 
\begin{eqnarray}
u_{1,\lambda}(x)&=&  \frac{\lambda\Vert W_2\Vert_2^{-1} - 
\sqrt{\lambda^2\Vert W_2\Vert_2^{-2} - 4b}}{2}\Vert W_2\Vert_2^{-1}W_2(x),
\\
u_{2,\lambda}(x)&=&  \frac{\lambda\Vert W_2\Vert_2^{-1} + 
\sqrt{\lambda^2\Vert W_2\Vert_2^{-2} - 4b}}{2}\Vert W_2\Vert_2^{-1}W_2(x).
\end{eqnarray}
}
We see from (1.13) and (1.14) that 
these two curves start from $(\lambda_0,\alpha_0)$ 
($\alpha_0:= \lambda_0\Vert W_2\Vert_2^{-2}\Vert W_2\Vert_\infty)$. Further, by Taylor expansion, we see that 
$\alpha_1(\lambda) = b\Vert W_2\Vert_2\lambda^{-1}(1 + o(1))$ and 
$\alpha_2(\lambda) = \Vert W_2\Vert_\infty
\Vert W_2\Vert_2^{-2}\lambda(1 + o(1))$ for $\lambda \gg 1$.

For the case $p > 1$ and $q = 1$ ($p \not= 2$),  it seems difficult to obtain such exact solutions 
$u_\lambda$ as (1.13)--(1.14). 
Therefore, we try to find the asymptotic shape of solutions $u_\lambda$ for $\lambda \gg 1$.  

\vspace{0.2cm}

\noindent
{\bf Theorem 1.3.} {\it Let $p > 1$ and $q = 1$. 
Put 
\begin{eqnarray}
\lambda_0:= (b(p-1))^{1/p}\frac{p}{p-1}\Vert W_p\Vert_p^{p-1}.
\end{eqnarray} 
(i) If $0 < \lambda < \lambda_0$, then there exists no solution of (1.1).

\noindent
(ii) If $\lambda = \lambda_0$, then there exists a unique solution of (1.1).

\noindent
(iii) If $\lambda > \lambda_0$, then there exist exactly two solutions $u_{1,\lambda}$ and 
$u_{2,\lambda}$ of (1.1). Moreover, for $\lambda \gg 1$, 
\begin{equation}
\left\{
\begin{array}{l}
\lambda:= \lambda_1(\alpha) = b\Vert W_p\Vert_\infty^{p-1}\alpha^{-(p-1)}
\left\{1 + b^{-1}\Vert W_p\Vert_p^p\Vert W_p\Vert_\infty^{-p}\alpha^p(1 + o(1))\right\},
\\
u_{1,\lambda}(x) = b^{1/(p-1)}\lambda^{-1/(p-1)}
\left(1 + 
\frac{1}{p-1}b^{1/(p-1)}\Vert W_p\Vert_p^p\lambda^{-p/(p-1)}
(1 + o(1))
\right)
W_p(x).
\end{array}
\right.
\end{equation}
\begin{equation}
\left\{
\begin{array}{l}
\lambda:= \lambda_2(\alpha) = \Vert W_p\Vert_p^p\Vert W_p\Vert_\infty^{-1}\alpha 
+ b\Vert W_p\Vert_\infty^{p-1}\alpha^{1-p} + o(\alpha^{1-p}),
\vspace{0.1cm}
\\
u_{2,\lambda}(x) =  \left\{\lambda \Vert W_p\Vert_p^{1-p} 
-b
\left(\lambda\Vert W_p\Vert_p^{1-p}\right)^{1-p}(1 + o(1))\right\}\Vert W_p\Vert_p^{-1}W_p(x).
\end{array}
\right.
\end{equation}
}

 \vspace{0.2cm}
 
 Finally, we treat (1.1) under the condition (1.2). 
 \vspace{0.2cm}
 
\noindent
{\bf Theorem 1.4.} {\it Assume (1.2) and $\lambda \gg 1$. 
Then there exist exactly two solutions $u_{1,\lambda}$ and 
$u_{2,\lambda}$ of (1.1). Moreover, for $\lambda \gg 1$, 
\begin{eqnarray}
u_{1,\lambda}(x) &=&  b^{q/(p-1)}\lambda^{-1/(p-1)}
\\
&&\times
\left\{1 + \frac{q}{p-1}b^{(pq-p+1)/(p-1)}\Vert W_p\Vert_p^p
\lambda^{-p/(p-1)}(1 + o(1))\right\}W_p(x),
\nonumber
\\
u_{2,\lambda}(x) &=& \left\{m\lambda^{1/(pq-p+1)} 
-\frac{bqm^{1-p}}{pq-p+1}\lambda^{(1-p)/(pq-p+1)}(1 + o(1))\right\}\Vert W_p\Vert_p^{-1}W_p(x).
\end{eqnarray}
where $m:= \Vert W_p\Vert_p^{(1-p)/(pq-p+1)}$.}

\vspace{0.2cm}

The remainder of this paper is organized as follows. In Section 2, we first introduce how to 
find the solutions of (1.1). Then we prove Theorem 1.1 by time map method. 
In Section 3, we mainly consider the existence of solutions of (1.1) under the condition (1.2). 
In Sections 4 and 5, we prove Theorems 1.3 and 1.4 by Taylor expansion and direct calculation.

 \section{Proof of Theorem 1.1}

In this section, let $b = 0$ in (1.1). 

\vspace{0.2cm}

\noindent
{\bf Lemma 2.1.} {\it For any $\lambda > 0$, (1.1) has a unique solution $u_\lambda$.}

\noindent
{\bf Proof.} We apply the argument [1, Theorem 2] to (1.1). 
 For a 
given $\lambda > 0$, we consider
\begin{equation}
\left\{
\begin{array}{l}
-w''(x)=\lambda w(x)^p, \enskip 
x \in I,
\vspace{0.1cm}
\\
w(x) > 0, \enskip x\in I.
\vspace{0.1cm}
\\
w(0) = w(1) = 0.
\end{array}
\right.
\end{equation}
Then it is clear that  
\begin{eqnarray}
w_\lambda(x):= \lambda^{1/(1-p)}W_p(x)
\end{eqnarray}
 is the unique solution of (2.1). 
For $t > 0$, we consider
\begin{eqnarray}
g(t):= t^{pq/2} - \Vert w_\lambda\Vert_p^{1-p}t^{(p-1)/2}.
\end{eqnarray} 
Then it is known from [1,Thorem 2] that if $g(t_\lambda) = 0$, then 
$u_\lambda:=\gamma w_\lambda$ ($\gamma:=  t_\lambda^{1/2}\Vert w_\lambda\Vert_p^{-1}$) 
satisfies (1.1). Indeed, by (2.1) and (2.3), we have
\begin{eqnarray}
&&-\left(\int_0^1 u_\lambda(x)^p dx \right)^q u_\lambda''(x)= 
-\Vert u_\lambda\Vert_p^{pq}u_\lambda''(x) 
\\
&&=- (\Vert \gamma w_\lambda\Vert_p)^{pq}\gamma w_\lambda''(x) = 
t_\lambda^{pq/2}\gamma\lambda w_\lambda(x)^p = (t_\lambda^{1/2}\Vert w_\lambda \Vert_p^{-1})^{p-1}
\gamma\lambda w_\lambda(x)^p
\nonumber
\\
&&= \lambda (\gamma w_\lambda(x))^p = \lambda u_\lambda(x)^p.
\nonumber
\end{eqnarray} 
On the other hand, assume that $u_\lambda$ satisfies (1.1).  Then we see that 
$u_\lambda = \Vert u_\lambda\Vert_p^{pq/(p-1)}w_\lambda$. Then we have 
$\Vert u_\lambda\Vert_p^{(p-1-pq)/(p-1)} = \Vert w_\lambda\Vert_p$. We put 
$t_\lambda ^{1/2}:= \Vert u_\lambda\Vert_p = \Vert w_\lambda\Vert_p^{(1-p)/(pq-p+1)}$. Then 
by (2.3), we have $g(t_\lambda) = 0$. Concequently, the solutions 
$\{t_{1,\lambda}, t_{2,\lambda}, \cdots, t_{k,\lambda}\}$ of (2.3) 
correspond to the solutions $\{u_{1,\lambda}, u_{2,\lambda}, \cdots, u_{k,\lambda}\}$. Therefore, 
if $g(t) = 0$ has a unique soution, then (1.1) also has a unique solution. By (1.2) and (2.3), we see that
\begin{eqnarray}
g(t) = t^{(p-1)/2}\left(t^{(pq-p+1)/2} - \Vert w_\lambda\Vert_p^{1-p}\right) = 0
\end{eqnarray}
has a unique positive solution $t_\lambda = \Vert w_\lambda\Vert_p^{2(1-p)/(pq-p+1)}$. Thus the proof is complete. 
\qed

\vspace{0.2cm}

\noindent
{\bf Proof of Theorem 1.1.} (i) Let $t_\lambda = \Vert w_\lambda\Vert_p^{2(1-p)/(pq-p+1)}$. 
By Lemma 2.1 and (2.2), we have 
\begin{eqnarray}
u_\lambda(x) = t_\lambda^{1/2}\Vert w_\lambda\Vert_p^{-1}w_\lambda(x) = t_\lambda^{1/2}
\Vert W_p \Vert_p^{-1}W_p(x).
\end{eqnarray}
By (2.2) and (2.5), we have 
\begin{eqnarray}
t_\lambda = \Vert w_\lambda\Vert_p^{2(1-p)/(pq-p+1)} = \lambda^{2/(pq-p+1)}
\Vert W_p\Vert_p^{2(1-p)/(pq-p+1)}.
\end{eqnarray}
 We apply the time map argument to (1.4). 
(cf. [10]). Since (1.4) is autonomous, we have 
\begin{eqnarray}
W_p(x) &=& W_p(1-x), \quad x \in [0, 1/2],
\\
W_p'(x) &>& 0, \quad x \in [0,1/2),
\\
\xi& :=& \Vert W_p\Vert_\infty = \max_{0\le x \le 1}W_p(x) = W_p(1/2).
\end{eqnarray}
By (1.4), for $0 \le x \le 1$, we have 
\begin{eqnarray}
\{W_p''(x) + W_p(x)^p\}W_p'(x) = 0.
\end{eqnarray}
By this and (2.10), we have 
\begin{eqnarray}
\frac12W_p'(x)^2 + \frac{1}{p+1}W_p(x)^{p+1} = \mbox{constant} = \frac{1}{p+1}W_p(1/2)^{p+1} 
= \frac{1}{p+1}\xi^{p+1}.
\end{eqnarray}
By this and (2.9), for $0 \le x \le 1/2$, we have 
\begin{eqnarray}
W_p'(x) = \sqrt{\frac{2}{p+1}(\xi^{p+1} - W_p(x)^{p+1})}.
\end{eqnarray}
By this, (2.8) and putting $\theta:= W_p(x)$, we have 
\begin{eqnarray}
\Vert W_p\Vert_p^p &=& 2\int_0^{1/2} W_p(x)^p dx 
= 2\int_0^{1/2} W_p(x)^p\frac{W_p'(x)}{\sqrt{\frac{2}{p+1}(\xi^{p+1} - W_p(x)^{p+1})}}dx
\\
&=& \sqrt{2(p+1)}\int_0^\xi \frac{\theta^p}{\sqrt{\xi^{p+1}-\theta^{p+1}}}d\theta
\qquad (\theta =\xi s)
\nonumber
\\
&=& \sqrt{2(p+1)}\xi^{(p+1)/2}\int_0^1 \frac{s^p}{\sqrt{1-s^{p+1}}}ds 
=  \sqrt{2(p+1)}\xi^{(p+1)/2}M_p.
\nonumber
\end{eqnarray}
By this, (2.6), and (2.7), we have 
\begin{eqnarray}
u_\lambda(x) = \lambda^{1/(pq-p+1)}\{(2(p+1))^{1/2}\xi^{(p+1)/2}M_p\}^{-q/(pq-p+1)}
W_p(x).
\end{eqnarray}
By putting $x = 1/2$ in (2.15), we have 
\begin{eqnarray}
\alpha = \lambda^{1/(pq-p+1)}\{(2(p+1))^{1/2}\xi^{(p+1)/2}M_p\}^{-q/(pq-p+1)}\xi.
\end{eqnarray}
By (2.13), we have 
\begin{eqnarray}
\frac12 &=& \int_0^{1/2} 1 dx = \int_0^{1/2} 
\frac{W_p'(x)}{\sqrt{\frac{2}{p+1}(\xi^{p+1} - W_p(x)^{p+1})}}dx 
\\
&=& \sqrt{\frac{p+1}{2}}\int_0^\xi \frac{1}{\sqrt{\xi^{p+1}-\theta^{p+1}}}d\theta
\qquad (\theta = \xi s)
\nonumber
\\
&=& \sqrt{\frac{p+1}{2}}\xi^{(1-p)/2}\int_0^1 \frac{1}{\sqrt{1-s^{p+1}}}ds 
= \sqrt{\frac{p+1}{2}}\xi^{(1-p)/2}L_p.
\nonumber
\end{eqnarray}
By this, we have 
\begin{eqnarray}
\xi = (2(p+1))^{1/(p-1)}L_p^{2/(p-1)}.
\end{eqnarray} 
By this, (2.15) and (2.16), we obtain (1.5) and (1.6). \qed

\noindent
(ii) Let $p = 1$. Then by (1.1), we have 
\begin{eqnarray}
-u_\lambda''(x) = \frac{\lambda}{\Vert u_\lambda\Vert_1^q}u_\lambda(x) = \pi^2u_\lambda(x).
\end{eqnarray}
It is clear that $u_\lambda(x) = \alpha\sin\pi x$, and 
\begin{eqnarray}
\Vert u_\lambda\Vert_1= \int_0^1 \alpha \sin\pi x dx = \frac{2\alpha}{\pi}.
\end{eqnarray}
By this and (2.19), we have 
\begin{eqnarray}
\lambda = \pi^2\Vert u_\lambda\Vert_1^q = 2^q\pi^{2-q}\alpha^q.
\end{eqnarray}
This implies (1.8). Thus the proof of Theorem 1.1 is complete. \qed

\vspace{0.2cm}

By (2.14) and (2.18), for $p > 1$, we obtain 

\begin{eqnarray}
\Vert W_p\Vert_p = (2(p+1))^{1/(p-1)}M_p^{1/p}L_p^{(p+1)/(p(p-1))}.
\end{eqnarray}
This implies (1.9). By (2.14) and (2.22), we obtain (1.10). \qed

\section{Proof of Theorem 1.2}

First, we consider (1.1) under the condition (1.2). 
The approach to find the solutions of (1.1) is a variant of  (2.3)--(2.4).
Namely,  we seek the solutions of (1.1)  of the form
\begin{eqnarray}
u_\lambda(x):= t \Vert w_\lambda\Vert_p^{-1} w_\lambda(x) = 
t \Vert W_p\Vert_p^{-1} W_p(x)
\end{eqnarray} 
for some $t > 0$. To do this, let $M(s):= (s + b)^q$. If we have solutions of (1.1) of the form (3.1), 
then since $\Vert u_\lambda\Vert_p = t$ by (3.1), we have 
\begin{eqnarray}
-M(\Vert u_\lambda\Vert_p^p)u_\lambda''(x) &=& -M(t^p)
\frac{t}{\Vert w_\lambda\Vert_p}w_\lambda''(x) 
\\
&=& -M(t^p)\frac{t}{\Vert w_\lambda\Vert_p}\lambda w_\lambda(x)^p
\nonumber
\\
&=& M(t^p)t^{1-p}\Vert w_\lambda\Vert_p^{p-1}\lambda u_\lambda(x)^p.
\nonumber
\end{eqnarray} 
By this, we look for $t$ satisfying
\begin{eqnarray}
 M(t^p)t^{1-p}\Vert w_\lambda\Vert_p^{p-1}=1.
\end{eqnarray}
Namely, we solve the equation 
\begin{eqnarray}
g(t):= (t^p + b)^q - \Vert w_\lambda\Vert_p^{1-p}t^{p-1} = 0. 
\end{eqnarray}
By this, we have 
\begin{eqnarray}
g'(t) &=& pq(t^p+b)^{q-1}t^{p-1} - (p-1)\Vert w_\lambda\Vert_p^{1-p}t^{p-2}
\\
&=& t^{p-2}\{pq(t^p+b)^{q-1}t - (p-1)\Vert w_\lambda\Vert_p^{1-p}\}=: t^{p-2}\tilde{g}(t).
\nonumber
\end{eqnarray}
By direct calculation, we see that $\tilde{g}(t)$ is 
strictly increasing for $t > 0$. Further, by (1.2) and (3.5), we have 
\begin{eqnarray*}
g'(t) \le pqt^{pq-1} - (p-1)\Vert w_\lambda\Vert_p^{1-p}t^{p-2} \to -\infty \quad 
(t \to 0). 
\end{eqnarray*}
Therefore, 
we see that there exists a unique $t_0 >0$ 
such that $g'(t) < 0$ for $0 < t_0 < t$, $g'(t_0) = 0$ and $g'(t) > 0$ for $t > t_0$. By using 
(3.4) and (3.5), we find that 
\begin{eqnarray}
g(t_0)&=& \frac{\Vert w_\lambda\Vert_p^{1-p}}{pqt_0}\left\{-(pq-p+1)t_0^p + b(p-1)\right\}.
\end{eqnarray}
If $g(t_0) < 0$, then there exists exactly $0 < t_1 < t_0 < t_2$ such that 
$g(t_1)=g(t_2) = 0$. If $g(t_0) > 0$, then (3.3) has no solutions. This idea will be used 
in the next sections. Proof of Theorem 1.2 is more simple, since $t_0$ is obtained explicitly. 

\vspace{0.2cm}

\noindent
{\bf Proof of Theorem 1.2.} Since $p=2, q=1$, by (3.4), we have 
\begin{eqnarray}
g(t) = t^2 + b - \lambda\Vert W_2\Vert_2^{-1}t = 0.
\end{eqnarray}
Then 
\begin{eqnarray}
t_{1,\lambda} &=& \frac{\lambda\Vert W_2\Vert_2^{-1} - \sqrt{\lambda^2\Vert W_2\Vert_2^{-2} - 4b}}
{2},
\\
t_{2,\lambda} &=& \frac{\lambda\Vert W_2\Vert_2^{-1} + \sqrt{\lambda^2\Vert W_2\Vert_2^{-2} - 4b}}{2}.
\end{eqnarray}
By this and (3.1), we obtain (i), (ii) and (iii). Thus the proof is complete. \qed

\section{Proof of Theorem 1.3}

In this section, let $p > 1, q=1$. 
Since $\Vert w_\lambda\Vert_p = \lambda^{-1/(p-1)} 
\Vert W_p\Vert_p$ by (2.2), we have from (3.5) that
\begin{eqnarray}
g'(t) = t^{p-2}\{pt - \lambda\Vert W_p\Vert^{1-p}(p-1)\}.
\end{eqnarray}  
We put 
\begin{eqnarray}
t_0 := \frac{p-1}{p}\lambda\Vert W_p\Vert_p^{1-p}.
\end{eqnarray}
Then $g'(t_0) = 0$. By this and (3.4), we have 
\begin{eqnarray}
g(t_0) &=& \left(\frac{p-1}{p}\lambda\Vert W_p\Vert_p^{1-p}\right)^p + b-
\lambda\Vert W_p\Vert_p^{1-p}\left(\frac{p-1}{p}\lambda\Vert W_p\Vert_p^{1-p}\right)^{p-1}
\\
&=& -\frac{1}{p-1}\left(\frac{p-1}{p}\lambda\Vert W_p\Vert_p^{1-p}\right)^p + b.
\nonumber
\end{eqnarray}
We put 
\begin{eqnarray}
\lambda_0:= (b(p-1))^{1/p}\frac{p}{p-1}\Vert W_p\Vert_p^{p-1}.
\end{eqnarray} 
By this, we see that (i)--(iii) are valid, since if $\lambda_0$ satisfies (4.4), then $g(t_0) = 0$ by (4.3). 
Further, $g(0) = b > 0$ and $g(t) > 0$ when $t \gg 1$. 

We now prove (1.17). We assume that $\lambda \gg 1$. 
Then there exists $0 < t_1 < t_0 < t_2$ which satisfy $g(t_1) = g(t_2) = 0$. Since $t_0 \to \infty$ 
as $\lambda \to \infty$, we see that $t_2 \to \infty$ as $\lambda \to \infty$. Then by  
(3.4), we have 
\begin{eqnarray}
t_2 = \lambda \Vert W_p\Vert_p^{1-p} + R,
\end{eqnarray}
where $R$ is the remainder term, and $R = o(\lambda)$. By (3.4), (4.5) and Taylor expansion, 
we have 
\begin{eqnarray}
&&\left(\lambda\Vert W_p\Vert_p^{1-p}\right)^p
\left(1 + \frac{R}{\lambda\Vert W_p\Vert_p^{1-p}}\right)^p + b - 
\left(\lambda\Vert W_p\Vert_p^{1-p}\right)^p
\left(1 + \frac{R}{\lambda\Vert W_p\Vert_p^{1-p}}\right)^{p-1} 
\\
&&= \left(\lambda\Vert W_p\Vert_p^{1-p}\right)^p
\left(1 + \frac{pR}{\lambda\Vert W_p\Vert_p^{1-p}}(1 + o(1))\right) + b 
\nonumber
\\
&&- \left(\lambda\Vert W_p\Vert_p^{1-p}\right)^p
\left(1 + \frac{R(p-1)}{\lambda\Vert W_p\Vert_p^{1-p}}(1 + o(1))\right) = 0.
\nonumber
\end{eqnarray} 
This implies that 
\begin{eqnarray}
R = -b
\left(\lambda\Vert W_p\Vert_p^{1-p}\right)^{1-p}(1 + o(1)).
\end{eqnarray}
By this, (3.1) and (4.5), for $\lambda \gg 1$, we have  
\begin{eqnarray}
u_{2,\lambda}(x) =  \left\{\lambda \Vert W_p\Vert_p^{1-p} 
-b
\left(\lambda\Vert W_p\Vert_p^{1-p}\right)^{1-p}(1 + o(1))\right\}\Vert W_p\Vert_p^{-1}W_p(x).
\end{eqnarray}
By putting $x = 1/2$ in (4.8), we have 
\begin{eqnarray}
\alpha = \left\{\lambda \Vert W_p\Vert_p^{1-p} 
-b
\left(\lambda\Vert W_p\Vert_p^{1-p}\right)^{1-p}(1 + o(1))\right\}\Vert W_p\Vert_p^{-1}
\Vert W_p\Vert_\infty.
\end{eqnarray}
By this, we obtain 
\begin{eqnarray}
\lambda &=& \Vert W_p\Vert_p^p\Vert W_p\Vert_\infty^{-1}\alpha 
+ b\Vert W_p\Vert_\infty^{p-1}\alpha^{1-p} + o(\alpha^{1-p}).
\end{eqnarray}
By this and (4.8), we obtain (1.17). 
We next prove (1.16). To do this, we 
consider the asymptotic behavior of $t_1$ as $\lambda \to \infty$. 
If there exists a constant $C > 0$ such that $C < t_1 < C^{-1}$. Then by (3.3), for 
$\lambda \gg 1$, we have 
\begin{eqnarray}
g(t_1) = t_1^p + b - \lambda\Vert W_p\Vert_p^{1-p}t_1^{p-1} < 0.
\end{eqnarray}
This is a contradiction, since $g(t_1) = 0$. If $t_1 \to \infty$ as $\lambda \to \infty$, then 
by (4.2) and (4.11), 
$t_1 > t_0$ for $\lambda \gg 1$. This is a contradiction. Therefore, $t_1 \to 0$ as 
$\lambda \to \infty$. By (3.4), we have 
\begin{eqnarray}
t_1^{p-1}(\lambda\Vert W_p\Vert_p^{1-p} - t_1) = b. 
\end{eqnarray}
By this and Taylor expansion, we have 
\begin{eqnarray}
t_1 &=& 
\left(\frac{b}{\lambda\Vert W_p\Vert_p^{1-p} - t_1}\right)^{1/(p-1)}
\\
&=& \Vert W_p\Vert_pb^{1/(p-1)}\lambda^{-1/(p-1)}
\left(\frac{1}{1-t_1\lambda^{-1}\Vert W_p\Vert_p^{p-1}}\right)^{1/(p-1)}
\nonumber
\\
&=& \Vert W_p\Vert_pb^{1/(p-1)}\lambda^{-1/(p-1)}
\left(1 + \frac{1}{p-1}\frac{t_1}{\lambda\Vert W_p\Vert_p^{1-p}}(1 + o(1))\right).
\nonumber
\end{eqnarray}
By this and (3.1), we have
\begin{eqnarray}
u_{1,\lambda}(x) &=& b^{1/(p-1)}\lambda^{-1/(p-1)}
\left(1 + 
\frac{1}{p-1}\frac{b^{1/(p-1)}\Vert W_p\Vert_p^p}{\lambda^{p/(p-1)}}
(1 + o(1))
\right)
W_p(x).
\end{eqnarray}
By this, we obtain 
\begin{eqnarray}
\alpha &=& b^{1/(p-1)}\lambda^{-1/(p-1)}
\left(1 + \frac{1}{p-1}\frac{b^{1/(p-1)}\Vert W_p\Vert_p^p}{\lambda^{p/(p-1)}}(1 + o(1))\right)\xi.
\end{eqnarray}
By this and (1.10), we have 
\begin{eqnarray}
\lambda_1 = b\Vert W_p\Vert_\infty^{p-1}\alpha^{-(p-1)}
\left\{1 + b^{-1}\Vert W_p\Vert_p^p\Vert W_p\Vert_\infty^{-p}\alpha^p(1 + o(1))\right\}.
\end{eqnarray}
By this and (4.14), we obtain (1.16). Thus the proof is complete. \qed

\section{Proof of Theorem 1.4}

In this section, we assume (1.2) and $\lambda \gg 1$. We put $k:= ((p-1)\Vert W_p\Vert_p^{1-p}/pq)^{1/(pq-p+1)}$. By (3.5), we have 
\begin{eqnarray}
t_0 = k\lambda^{1/(pq-p+1)}(1 + o(1)). 
\end{eqnarray}
By this, (1.2) and (3.4), we see that $g(t_0) < 0$. Then there exists  $0 < t_1 < t_0 < t_2$ 
such that $g(t_1) = g(t_2) = 0$. By (5.1), we see that $t_2 \to \infty$ as $\lambda \to \infty$.
We first prove (1.19).  We recall that $m:= \Vert W_p\Vert_p^{(1-p)/(pq-p+1)}$. By (3.4), we have 
\begin{eqnarray}
t_2 = m\lambda^{1/(pq-p+1)} + r,
\end{eqnarray}
 where $r$ is the remainder term satisfying $r = o(\lambda^{1/(pq-p+1)})$. 
It is clear that  (3.4) is equivalent to 
 \begin{eqnarray}
 (t_2^p + b)^q = \lambda\Vert W_p\Vert_p^{1-p}t_2^{p-1}.
 \end{eqnarray}
 By this, (5.2) and 
 Taylor expansion, we have 
\begin{eqnarray}
 \mbox{r.h.s. of (5.3)} &=& \lambda\Vert W_p\Vert_p^{1-p}
 \left(m\lambda^{1/(pq-p+1)} + r\right)^{p-1} 
 \\
 &=& \lambda\Vert W_p\Vert_p^{1-p}m^{p-1}\lambda^{(p-1)/(pq-p+1)}
 \left(1 + (p-1)\frac{r}{m \lambda^{1/(pq-p+1)}}(1 + o(1)\right).
 \nonumber
\end{eqnarray}
By (5.2) and Taylor expansion, we have 
\begin{eqnarray}
\mbox{l.h.s. of (5.3)} &=& t_2^{pq}\left(1 + \frac{b}{t_2^{p}}\right)^q 
= t_2^{pq}\left(1 + \frac{bq}{t_2^{p}}(1 + o(1))\right)
\\
&=& \left(m\lambda^{1/(pq-p+1)} + r\right)^{pq}
\left\{1 + bq(m\lambda^{1/(pq-p+1)} + r)^{-p}(1+o(1))\right\}
\nonumber
\\
&=& (m\lambda^{1/(pq-p+1)})^{pq}
\left(1 +\frac{r}{m\lambda^{1/(pq-p+1)}}(1 + o(1))\right)^{pq}
\nonumber
\\
&&\times
\left\{
1 + bq(m\lambda^{1/(pq-p+1)})^{-p}(1 + o(1))
\right\}.
\nonumber
\end{eqnarray}
By the definition of $m$, we see that the leading terms of (5.4) and (5.5) coinside each other. 
By (5.4) and (5.5), we have 
\begin{eqnarray}
pq\frac{r}{m\lambda^{1/(pq-p+1)}} + \frac{bq}{m^p\lambda^{p/(pq-p+1)}} 
= (p-1)\frac{r}{m\lambda^{1/(pq-p+1)}}.
\end{eqnarray}
This implies that 
\begin{eqnarray}
r = -\frac{bqm^{1-p}}{pq-p+1}\lambda^{(1-p)/(pq-p+1)}(1 + o(1)). 
\end{eqnarray}
By this and (5.2), for $\lambda \gg 1$, we have 
\begin{eqnarray}
t_2 = \left\{m\lambda^{1/(pq-p+1)} 
-\frac{bqm^{1-p}}{pq-p+1}\lambda^{(1-p)/(pq-p+1)}(1 + o(1))\right\}.
\end{eqnarray}
By this and (3.1), we have 
\begin{eqnarray}
u_{2,\lambda}(x) = \left\{m\lambda^{1/(pq-p+1)} 
-\frac{bqm^{1-p}}{pq-p+1}\lambda^{(1-p)/(pq-p+1)}(1 + o(1))\right\}\Vert W_p\Vert_p^{-1}W_p(x).
\end{eqnarray}
This implies (1.19). We next show (1.18). We consider the asymptotic behavior of $t_1$ as $\lambda \to \infty$. By the same argument as that in Section 4, we find that $t_1 \to 0$ as $\lambda \to \infty$. By 
(5.3), we have 
\begin{eqnarray}
\lambda\Vert W_p\Vert_p^{1-p}t_1^{p-1} = b^q(1 + o(1)).
\end{eqnarray}
This implies that 
\begin{eqnarray}
t_1 = b^{q/(p-1)}\Vert W_p\Vert_p\lambda^{-1/(p-1)}(1 + \eta),
\end{eqnarray}
where $\eta$ is the remainder term. Then by Taylor expansion and (5.11), we have 
\begin{eqnarray}
\mbox{l.h.s. of (5.3)} &=& (b + t_1^p)^q = b^q(1 + b^{-1}t_1^p)^q = b^q(1 + b^{-1}qt_1^p + o(t^p)),
\\
\mbox{r.h.s. of (5.3)} &=& b^q(1 + (p-1)\eta + o(\eta)).
\nonumber
\end{eqnarray}
By this, we have 
\begin{eqnarray}
\eta &=& \frac{q}{p-1}b^{-1}t_1^p = \frac{q}{p-1}b^{(pq-p+1)/(p-1)}\Vert W_p\Vert_p^p
\lambda^{-p/(p-1)}(1 + o(1)).
\end{eqnarray}
By this and (5.11), we have 
\begin{eqnarray}
u_{1,\lambda}(x) &=&  b^{q/(p-1)}\lambda^{-1/(p-1)}
\\
&&\times
\left\{1 + \frac{q}{p-1}b^{(pq-p+1)/(p-1)}\Vert W_p\Vert_p^p
\lambda^{-p/(p-1)}(1 + o(1))\right\}W_p(x).
\nonumber
\end{eqnarray}
This implies (1.18). Thus the proof is complete. \qed


\begin{thebibliography}{20}
\labelsep=1em\relax

\bibitem[1]{1} Alves, C. O.; Corr\'ea, F. J. S. A.; Ma, T. F. Positive solutions for a quasilinear elliptic equation of Kirchhoff type. Comput. Math. Appl. 49 (2005), 85--93.
\bibitem[2]{2}  Arcoya, David; Leonori, Tommaso; Primo, Ana, Existence of solutions for semilinear nonlocal elliptic problems via a Bolzano theorem. Acta Appl. Math. 127 (2013), 87--104.
\bibitem[3]{3} F. J. S. A. Corr\^{e}a, On positive solutions of nonlocal and nonvariational elliptic problems, Nonlinear Anal. 59 (2004), 1147--1155.
\bibitem[4]{4}Corr\^{e}a, Francisco Julio S. A.; de Morais Filho, Daniel C. On a class of nonlocal elliptic problems via Galerkin method. J. Math. Anal. Appl. 310 (2005), no. 1, 177--187.
\bibitem[5]{5} {J. M. Fraile, J. L\'opez-G\'omez and J. Sabina de Lis}, 
{On the global structure of the set of positive 
solutions of some semilinear elliptic boundary value 
problems}, {J. Differential Equations} {\bf 123} (1995), 180--212. 
\bibitem[6]{6} B. Gidas, W. M. Ni and L. Nirenberg, Symmetry and related properties via the maximum principle. Comm. Math. Phys. 68 (1979), 209--243.
\bibitem[7]{7} C.S. Goodrich, A topological approach to nonlocal elliptic partial differential equations on an annulus. Math. Nachr. 294 (2021), 286--309.
\bibitem[8]{8} {A.A. Lacey}, Thermal runaway in a non-local problem modelling Ohmic heating. I. Model derivation and some special cases. European J. Appl. Math. 6 (1995), 127--144.
\bibitem[9]{9} {A.A. Lacey,} Thermal runaway in a non-local problem modelling Ohmic heating. II. General proof of blow-up and asymptotics of runaway. European J. Appl. Math. 6 (1995), 201--224.
\bibitem[10]{10} {T. Laetsch}, {
The number of solutions of a nonlinear 
two point boundary value problem}, 
Indiana Univ. Math. J. {\bf 20}  1970/1971 1--13.
\bibitem[11]{11} 
R. Sta\'nczy, Nonlocal elliptic equations, Nonlinear Anal. 47 (2001), 3579--3584.



\end{thebibliography}
\end{document}